\documentclass[11pt]{amsart}
\usepackage{amssymb}
\usepackage{verbatim}

\newtheorem{theorem}{Theorem}[section]
\newtheorem{lemma}[theorem]{Lemma}
\newtheorem{proposition}[theorem]{Proposition}
\newtheorem{corollary}[theorem]{Corollary}
\newtheorem{remark}[theorem]{Remark}

\theoremstyle{definition}

\numberwithin{equation}{section}

\begin{document}

\title[Some split 
symbol algebras of prime degree]{Some split 
symbol algebras of prime degree}

\author{Diana Savin}
\address{Faculty of Mathematics and Computer Science, Ovidius University\\
Bd. Mamaia 124, 900527, Constanta, Romania}
\email{savin.diana@univ-ovidius.ro; dianet72@yahoo.com}

\author{Vincenzo Acciaro}
\address{Dipartimento di Economia\\ Universit\`a di Chieti--Pescara\\
Viale Pindaro 42\\  65127 Pescara, Italy}
\email{v.acciaro@unich.it}


\date{}


\begin{abstract}
Let $p$ be an odd prime,  let  $K=\mathbb{Q}(\epsilon)$ where $\epsilon$ is a primitive cubic root  of unity,  and let $L$ be the Kummer field $\mathbb{Q}\left(\epsilon, \sqrt[3]{\alpha}\right)$. In this paper we obtain a characterization of the splitting behavior of the symbol algebras   
 $\left( \frac{\alpha ,p}{K,\epsilon }\right)$
 and 
$\left( \frac{\alpha ,p^{h_{p}}}{K,\epsilon}\right)$, where $h_{p}$ is the order in the class group $Cl\left(L\right)$  of a prime ideal of $\mathcal{O}_L$ which divides $p\mathcal{O}_L.$ 
\end{abstract}

\maketitle

 \section{introduction}
Let $n\geq 3$ be an arbitrary positive integer,  and let $F$ be a field with $char(F) \nmid n$.   Let 
 $\xi$ be a primitive $n$-th root of unity in $F$. 
 If   $a,b$ $\in F\backslash \{0\}$, 
 the algebra $A$ over $F$ generated by two elements $x$ and $y$ satisfying%
\begin{equation*}
x^{n}=a,y^{n}=b,yx=\xi xy
\end{equation*}
is called a  {symbol algebra} and it is denoted by $\left( \frac{a,~b}{%
F,\xi }\right)$.
When $n = 2$ we obtain the well known generalized quaternion algebra over the field $F$, and indeed a symbol algebra is a natural generalization of a quaternion algebra.
Quaternion algebras and symbol algebras are central simple algebras of dimension $n^{2}$ over the base  field $F$.
The results about  quaternion algebras and symbol algebras have strong connections with number theory (especially with the ramification theory in algebraic number fields).
Different criteria are known for a quaternion algebra or a symbol algebra to split 
\cite{GiSz06, Lam04,Vo10}.
Explicit conditions for a quaternion algebra over the field of rationals numbers to be split or else a division algebra were studied in \cite{AlBa04}. In the paper \cite{Sa16} we investigated  the splitting behavior of quaternion algebras over  quadratic fields, and of specific  symbol algebras over cyclotomic fields. In   \cite{Sa17}  we found a sufficient condition for a quaternion algebra over a quadratic field to split.
Next, in  \cite{AcSaTaZe19}  we gave necessary and sufficient conditions for a quaternion algebra   $H(\alpha,m)$ to split    over a quadratic field $K$, and in   
\cite{AcSaTaZe191}  we obtained
a complete characterization of division quaternion algebras $H(p,q)$,
where $p,q$ are prime integers, 
over the
composite $K$ of $n$ quadratic number fields.

In this paper we study some symbol algebras of prime degree over very specific cyclotomic fields.
In  Section 2 we recall some useful results about symbol algebras, cyclotomic fields and Kummer fields which we will use later. In Section 3 we find conditions for a symbol algebra of prime degree over a cyclotomic field to split. The main result of this article is:\\
\smallskip\\
\textbf{Main Theorem (Theorem 3.7)} 
	\textit{Let} $\epsilon $  \textit{be a primitive cubic
		root }    \textit{of unity and let }$K=\mathbb{Q}\left( \epsilon
	\right) $. \textit{Let} $\alpha \in
	K^{\ast },$ \textit{let} $p$  \textit{be a prime rational integer,} $p\neq 3$ 
	\textit{and let } $L=K\left( \sqrt[3]{\alpha }\right) $  \textit{be the Kummer field}.  {Let} $\mathcal{O}_L$  \textit{be the ring of integers of the field} $L.$
	\textit{Let} $Cl\left(L\right)$  \textit{be the ideal class group of the ring} $\mathcal{O}_L,$  \textit{let } $h_{L}$  \textit{be the class number of} $L$  \textit{and let} $h_{p}$  \textit{be the order of the class of a prime ideal  $P$ in} $\mathcal{O}_L$, \textit{which divides} $p\mathcal{O}_L,$  \textit{in the group} $Cl\left(L\right).$  \textit{Then, there exists a unit} $u$ $\in$ $U\left(\mathbb{Z}\left[ \epsilon %
	\right]\right)$
	\textit{such that the symbol algebra} $A=\left( \frac{\alpha, u\cdot p^{h_{p}}}{K,\epsilon }\right) $
	\textit{splits if and only if} $\alpha $  \textit{is a cubic residue modulo} $p.$\\
\smallskip\\
 Our results have been computationally validated
with the aid of the computer algebra package MAGMA.
\section{Some basic results}
We   recall here the decomposition behavior  of a prime integer
 in the ring of integers of a cyclotomic field
 and in the ring of integers  of a Kummer field. We will use these facts later to prove our results.
\begin{proposition}
\label{prop21} \cite{IrRo92}  {Let } $l\geq 3$ {\ be an integer,}
$\xi $ {be a primitive root of unity of order   \thinspace }$l.$  {Let }$K=\mathbb{Q}\left(\xi\right)$
 {and let} $\mathcal{O}_K$  {be the ring of integers of the cyclotomic field} $K.$  {Then} $\mathcal{O}_K=
\mathbb{Z}\left[\xi\right].$
\end{proposition}
\begin{theorem}
\label{teo22} \cite{IrRo92} {Let }$l\geq 3$ {be an integer,}%
  {\ and let  }$\xi $ {be a primitive root of unity
of order   \thinspace }$l$. If $p$ {is a prime
  number which does not divide }$l$  and $f$%
 {\ is the smallest positive integer such that }$p^{f}\equiv 1$ {mod }$l,$ {then we have }$p\mathbb{Z}[\xi ]=P_{1}P_{2}....P_{r},$ 
 {where }$r= {\varphi \left( l\right) }/{f},$ where $\varphi $  is
the Euler's function and $P_{j},\,j=1,...,r$ {\ are different prime
ideals in the ring }$\mathbb{Z}[\xi].$
\end{theorem}
\begin{theorem}
\label{teo23} \cite{Lem00} {Let }$\xi $ {\ be a primitive
root of unity of order  }$l,$ {where }$l$ {\ is a prime
 number, and let} $\mathcal{O}_L$   {be the ring of integers of the Kummer
field} $L=\mathbb{Q}(\xi ,\sqrt[l]{\mu}),$ where $\mu$$\in$$\mathbb{Q}$. {A prime ideal} $P$ {of} $\mathbb{Z}[\xi ]$ {\ decomposes in} $\mathcal{O}_L$ as follows: 
\begin{itemize}
\item
 it is equal to the $l${-power of a prime ideal of} $\mathcal{O}_L,$  {if the} $l${-power character }   $\left( \frac{\mu }{P}%
\right) _{l}=0;$
 \item
 it is a prime ideal of $\mathcal{O}_L,$ {if} $\left( \frac{\mu }{P}%
\right) _{l}$ {is a root of order }$l$ {of unity, different
from }$1$;
\item
 it is equal to the  product of $l$ {different prime ideals of} $\mathcal{O}_L,$ {if} $\left( \frac{\mu }{P}\right) _{l}=1.$
\end{itemize}
\end{theorem}
We recall now the definition of the $l${-power character} $\left( \frac{\mu }{P} 
\right) _{l}.$ If $\mu$$\not\in$$P,$ there is a unique integer $c$ (modulo $l$) such that
$\mu ^{\frac{N\left(P\right)-1}{l}}$$\equiv$ $\xi^{c}$ (mod $P$).\\
The $l${-power character} $\left( \frac{\mu }{P} \right) _{l}$ is defined as:\\
i) if $\mu \in P$  then $\left( \frac{\mu }{P} \right) _{l}$$=0$;\\
ii) if $\mu \not\in P$ then  $\left( \frac{\mu }{P} \right) _{l}$ is the unique $l$th root of unity such that 
$$\mu ^{\frac{N\left(P\right)-1}{l}}\equiv \left( \frac{\mu }{P} \right) _{l} (mod\; P).$$
Moreover $\left( \frac{\mu }{P} \right) _{l}=1$ if and only if $\mu \not\in P$ and the congruence $x^{l}\equiv \mu$
(mod $P$) is solvable in $\mathbb{Z}[\xi ].$\\
\smallskip\\ 
Let's recall now some results about central simple algebras.
Let $A$ be a central simple algebra over a field $K.$
Then the dimension $n$ of $A$\ over $K$ is a square; its positive square root 
  is called the degree of the algebra $A$. 
If the equations $ax=b,\,ya=b$ have unique solutions for all $a,b \in A$, with $a \neq 0$, 
then the algebra $A$ is called  a division algebra. $A$ is a division algebra if and only if $A$
has no  zero divisors ($x\neq 0,y\neq 0\Rightarrow xy\neq 0$).
Let $K\subseteq L$ be a field extension and let $A$ be a central simple algebra over $K.$ Then:
\begin{itemize}
\item
$A$ is called  {split} by $K$ if $A$ is isomorphic to a full
matrix algebra over $K$;
\item
$A$ is called  {split} by $L$, and $L$ is called
a  {splitting field} for $A$, if $A\otimes _{K}L$
is a full matrix algebra over $L.$
\end{itemize}
The following splitting criteria for symbol algebras is known:
\begin{theorem}
\label{teo24} \cite{GiSz06} {Let }$K$ {be a field which contains a primitive n-th root of unity }%
$\xi$,
 and let  $a
, b \in K^{\ast }. $ {\ Then the following statements are equivalent:%
}
\begin{itemize}
\item
 {The symbol algebra }$A=\left( \frac{a, b}{K,\xi }%
\right) $ {\ is split.}
\item
 { The element }$b$ {is a norm from the extension }$%
K\subseteq K(\sqrt[n]{a}). $ 
\end{itemize}
\end{theorem}
For symbol algebras of prime degree the following is true:
\begin{remark}
\label{rem25} \cite{Le05}  {Let} $n\geq3$  {be a positive integer}, {and let} $\xi $  {be a
primitive} $n${-th root of unity. Let} $K$  {be a field of {characteristic} $\neq 2$  which contains} $\xi$ and let
$\alpha,\beta\in$$K^{*}$.  {If} $n$  {is prime, then the symbol algebra} $\left( \frac{\alpha,~\beta}{%
K,\xi }\right)$  {is either split or a division algebra.}
\end{remark}
\begin{lemma}
\label{lemma26} \cite{Dr}  {Let} $n\geq 2$  {be a positive integer,}   {and let} $\xi $  {be a primitive root of  unity of order} $n.$  
{Let} $K$  {be a field  which contains} $\xi$.
 {Let} $Br\left(K\right)$  {be the Brauer group of the field} $K$  {and let} $_{n}Br\left(K\right)$ {be the} $n$-{torsion component of} $Br\left(K\right).$ {Then, the assignment} $\left(\alpha,~\beta\right)\longmapsto$$\left( \frac{\alpha,~\beta}{%
K,\xi }\right)$  {induces a} $\mathbb{Z}-${bilinear map} 
$$K^{*}/\left(K^{*}\right)^{n} \times K^{*}/\left(K^{*}\right)^{n}\mapsto _{n}Br\left(K\right).$$
\end{lemma}
\section{Symbol algebras which split over specific cyclotomic fields}
In this section we study the symbol algebras $\left( \frac{\alpha ,p^{c}}{K,\xi }\right)$ 
when  $\xi $ is a primitive root of  unity of prime order $q$ and $L$ is the Kummer field $L=\mathbb{Q}\left(\xi, \sqrt[q]{\alpha} \right)$, for some  particular values of 
 $c$. 
We start with a small remark about such algebras:
\begin{remark}
\label{rem31}  {Let} $n\geq 3$  {be a positive integer,}  {let} $p$ {be a prime positive integer} {and let} $K=\mathbb{Q}\left(\xi\right)$
where $\xi$  {is a primitive
root of unity of order} $n$.   {Let} $\alpha \in
K^{\ast }$  {and let $L$ be the Kummer field} $K\left( \sqrt[n]{\alpha }\right).$
 {Then, the symbol algebras} $A=\left( \frac{\alpha ,p^{n}}{K, \xi}\right)$%
splits.
\end{remark}
\begin{proof}   
By Lemma 2.6  the symbol algebra $A=\left( \frac{\alpha ,p^{n}}{K, \xi}\right)$
lies in the same class of the algebra
$ \left( \frac{\alpha ,1 }{K, \xi}\right)$
in the Brauer group of $K$.
But $ \left( \frac{\alpha ,1 }{K, \xi}\right)$   splits by Theorem 2.4,
since 1 is always a norm.
\end{proof}
In the paper \cite{FlSa15} the authors obtained some results about the symbol algebras of the form $\left( \frac{%
\alpha ,p^{h_{L}}}{K,\xi }\right)$. In that paper there is a small
mistake, which we fix in the next proposition and in its corollary.
\begin{proposition}
\label{prop32}  \cite[Prop. 4.1]{FlSa15}
 {Let} $\epsilon $  {be a primitive cubic
root  of unity and let }$K=\mathbb{Q}\left( \epsilon
\right)$.   {Let} $\alpha \in
K^{\ast }$ be a cubic residue modulo $p$, where $p\neq 3$    {is a prime   integer}.
 {Let} $h_{L}$  {be the class number of
the Kummer field } $L = K\left( \sqrt[3]{\alpha }\right) .$  {Then, there exists a unit} $u$ $\in$ $U\left(\mathbb{Z}\left[ \epsilon\right]\right)$  { such that the symbol algebra} $A=\left( \frac{\alpha, u\cdot p^{h_{L}}}{K,\epsilon }\right) $%
 {\ splits.}
\end{proposition}
We recall that when $\left(R,+,\cdot\right)$ is a unitary ring with unity   $1$, then the set of the invertible elements of $\left(R,+,\cdot\right)$ is denoted by $U\left(R\right)$, i.e. $U\left(R\right)=\left\lbrace x\in R\;| \;\exists x^{'}\in R\; such \; that\; x\cdot x^{'}=x^{'}\cdot x=1\right\rbrace. $\\
\begin{corollary}
\label{cor33}
 \cite[Cor. 4.2]{FlSa15}  {Let} $q$  {be an odd prime 
integer, let } $\xi $  {be a primitive root of  unity of order} $q$,  {and let } $K=\mathbb{Q}\left( \xi \right) $.
  {Let}  $p\neq q$  {be a prime 
integer and let } $\alpha \in K^{\ast }$
 be a $q$  {power residue modulo} $p.$  {Let} $h_{L}$  {be
the class number of the Kummer field
} $L = K\left( \sqrt[q]{\alpha }\right)$.  {Then, there exists a unit} $u$ $\in$ $U\left(\mathbb{Z}\left[\xi %
\right]\right)$
 such that the symbol algebra $A=\left( \frac{\alpha, u\cdot p^{h_{L}}}{K, \xi}\right) $%
 {\ splits.}\\  
\end{corollary}
In the next two results we  show how to generalize the previous results, 
by replacing  $h_L$ with some divisors $h$ of  $h_L$.  
 \begin{proposition}
\label{prop35}  {Let} $q\geq 3$  {be a prime positive
integer},    {and} $K=\mathbb{Q}\left(\xi\right)$ where
 $\xi $  {is a primitive root of order} $q$  {of
unity.}  {Let} $\alpha \in
K^{\ast },$ $p$  {be}  {a prime rational integer,} $p\neq q$ 
 {and let } $L=K\left( \sqrt[q]{\alpha }\right) $  {be}  {%
the Kummer field such that} $\alpha $  {is a} $q$-th  {power residue modulo} $p.$  {Let} $\mathcal{O}_L$  {be the ring of integers of the field} $L.$
 {Let} $Cl\left(L\right)$  {be the ideal class group of the ring} $\mathcal{O}_L,$  {let }$h_{L}$  {be the class number of} $L$  {and let}
$h_{p}$ be the
order in the 
class group $Cl\left(L\right)$  of the class of a prime ideal in $\mathcal{O}_L$ that
divides $p\mathcal{O}_L.$ 
{Then}:
\begin{enumerate}
\item
{there exists a unit} $u$ $\in$ $U\left(\mathbb{Z}\left[\xi %
\right]\right)$
 { such that the symbol algebra} $A=\left( \frac{\alpha, u\cdot p^{h_{p}}}{K,\xi}\right) $%
 {\ splits;}
 \item there exists a unit $u$ $\in$ $U\left(\mathbb{Z}\left[\xi %
 \right]\right)$
 such that the symbol algebra $A=\left( \frac{\alpha, u\cdot p^{gcd\left( h_{p},q\right) }}{K,\xi}\right) $
 splits.
 \end{enumerate}
\end{proposition}
A particular case of Proposition \ref{prop35} is the following corollary.
\begin{corollary}
\label{cor34} 
 {Let} $\epsilon $  {be a primitive cubic
root  of unity and let }$K=\mathbb{Q}\left( \epsilon
\right)$.   {Let} $\alpha \in
K^{\ast }$  {be a cubic residue modulo} $p$  {with} $p\neq 3$  {a prime integer}.
 {Let} $h_{L}$  {be the class number of
the Kummer field } $L = K\left( \sqrt[3]{\alpha }\right) .$
Let $h_{p}$ be the
order in the 
class group $Cl\left(L\right)$  of the class of a prime ideal in $\mathcal{O}_L$ that
divides $p\mathcal{O}_L.$ 
{Then}:
\begin{enumerate}
\item
{there exists a unit} $u$ $\in$ $U\left(\mathbb{Z}\left[ \epsilon %
\right]\right)$
 { such that the symbol algebra} $A=\left( \frac{\alpha, u\cdot p^{h_{p}}}{K,\epsilon }\right) $ {splits;}
\item {there exists a unit} $u$ $\in$ $U\left(\mathbb{Z}\left[\epsilon %
\right]\right)$
{ such that the symbol algebra} $A=\left( \frac{\alpha, u\cdot p^{gcd\left( h_{p},3\right) }}{K,\epsilon }\right) $%
{\ splits.}
\end{enumerate}
 \end{corollary}
Since Proposition \ref{prop35} is 
more general than Corollary \ref{cor34}, 
we will prove only Proposition \ref{prop35}.\\ 
\smallskip\\
Before going into the proof of Proposition \ref{prop35} we would like to show   some examples of split symbol
algebras, which satisfy the hypotheses of Corollary \ref{cor34}, which were produced by using the computer
algebra package MAGMA \cite{Mag}:
\begin{itemize}
\item
Let  $K=\mathbb{Q}\left(\epsilon\right) $, where $\epsilon^{3}=1, \epsilon\neq 1$; the class number of Kummer field $L=\mathbb{Q}\left(\epsilon,\sqrt[3]{43 }\right) $ is $48$. Let $p=23.$
The  ideal $23\mathcal{O}_L$   decomposes into the product of three prime ideals of $\mathcal{O}_L$, \ $ P_{1}, P_{2}   , P_{3}     $ in the notations of our example. We denote with
$\left[I\right]$ be the class of the ideal $I$ in the class group of $L.$ We have $h_{p}=ord\left(\left[P_{1}\right]\right)=
ord\left(\left[P_{2}\right]\right)=ord\left(\left[P_{3}\right]\right)=12$ and $gcd\left(h_{p},3\right)=3.$ The norm equation $23^{12}=N_{L/\mathbb{Q}\left(\epsilon \right)}%
\left( a\right) $ has solutions. Also the norm equation $23^{3}=N_{L/\mathbb{Q}\left(\epsilon \right)}%
\left( a\right) $ has solutions.  
But the norm equations $23^{2}=N_{L/\mathbb{Q}\left(\epsilon \right)}%
\left( a\right), $ $23=N_{L/\mathbb{Q}\left(\epsilon \right)}\left( a\right)$  do not have any. 
This example agrees with the assertion of Corollary \ref{cor34}.
\item
Let $K=\mathbb{Q}\left( \epsilon \right) $, where $\epsilon^{3}=1, \epsilon\neq 1$; the class number of Kummer field $L=\mathbb{Q}\left(\epsilon,\sqrt[3]{43 }\right) $ is $48$. Let $p=11.$ The ideal $11\mathcal{O}_L$   decomposes into the product of three prime ideals of $\mathcal{O}_L$, \ $P_{1} , P_{2}, P_{3}$.
 We have $h_{p}=ord\left(\left[P_{1}\right]\right)=
ord\left(\left[P_{2}\right]\right)=ord\left(\left[P_{3}\right]\right)=2$ and $g.c.d\left(h_{p},3\right)=1.$
 The norm equation $11^{2}=N_{L/\mathbb{Q}\left(\epsilon \right)}%
\left( a\right) $ has solutions and also the norm equation $11=N_{L/\mathbb{Q}\left(\epsilon \right)}%
\left( a\right) $ has solutions.
This example agrees again with the assertion of Corollary \ref{cor34}.
\item
Let  $K=\mathbb{Q}\left(\epsilon\right) $, where $\epsilon^{3}=1, \epsilon\neq 1$; the class number of Kummer field $L=\mathbb{Q}\left(\epsilon,\sqrt[3]{11 }\right) $ is $4$. Let $p=19\equiv 1$ (mod $3$).
The  ideal $19\mathcal{O}_L$   decomposes into the product of six prime ideals of $\mathcal{O}_L$, $ P_{11}, P_{12}, P_{13}, P_{21}, P_{22}, P_{23}$ in the notations of our example. We denote with
$\left[I\right]$ be the class of the ideal $I$ in the class group of $L.$ We have $h_{p}=ord\left(\left[P_{11}\right]\right)=
ord\left(\left[P_{12}\right]\right)=ord\left(\left[P_{13}\right]\right)=ord\left(\left[P_{21}\right]\right)=
ord\left(\left[P_{22}\right]\right)=ord\left(\left[P_{23}\right]\right)=2$ and $g.c.d\left(h_{p},3\right)=1.$ The norm equation $19^{2}=N_{L/\mathbb{Q}\left(\epsilon \right)}%
\left( a\right) $ has solutions, and also the norm equation $19=N_{L/\mathbb{Q}\left(\epsilon \right)}%
\left( a\right)$  has solutions. 
This example agrees with the assertion of Corollary \ref{cor34}.
\end{itemize}
We are now going to prove Proposition \ref{prop35}:
\begin{proof}  
\begin{enumerate}
\item
	It is known that $\mathcal{O}_K=\mathbb{Z}\left[\xi %
	\right]$ and it is a Dedekind ring, therefore, any nonzero ideal of $\mathcal{O}_K$ decomposes uniquely into a product of prime ideals of the ring $\mathcal{O}_K.$
	We split the proof into two cases:
	\begin{itemize}
		\item
$p$ is a primitive root modulo $q$. \\
From Theorem \ref{teo22} it follows that $p$ remains prime in the ring $\mathbb{Z}\left[ \xi\right].$
		Since  the $q$-th power character   $\left( \frac{\alpha }{p}%
		\right) _{q}$ is equal to   $1$, from  Theorem 2.3 it follows that
		we have the following decomposition of the ideal $p\mathcal{O}_{L}$
		as a product of prime ideals  in $\mathcal{O}_{L}$:
		$$p\mathcal{O}_{L}=P_{1}\cdot P_{2}\cdot...\cdot P_{q},$$ 
		Let $h_{p}$ be the order of the ideal class of $P_{1}$ in the group $Cl\left(L\right).$ Since the ideals $P_{1},$ $P_{2},$..., $P_{q}$ are  conjugate under the action of the Galois group Gal($L/K$), it follows that the ideal classes of  $P_{1},$ $P_{2},$..., $P_{q}$ have the same order $h_{p}$ in the group $Cl\left(L\right).$
		Now $\left( p\mathcal{O}%
		_{L}\right) ^{h_{p}}=P_{1}^{h_{p}}P_{2}^{h_{p}}\cdot...\cdot P_{q}^{h_{p}}=(\beta\mathcal{O}_{L})   
		(\sigma(\beta) \mathcal{O}_{L})\cdot...\cdot (\sigma^{q-1}(\beta)\mathcal{O}_{L})$ 
		for some $\beta  \in  \mathcal{O}_{L} $,
		where $\sigma$ is a generator of $Gal(L/K)$.
		Hence
		$p^{h_{p}}    \mathcal{O}_{L }= N_{L/K}%
		\left( \beta\right)  \mathcal{O}_{L}.$ This means that $p^{h_{p}}  $   and    $N_{L/K}(\beta)$  differ by a unit of $\mathcal{O}_{L}$,
		but, since  $p^{h_{p}}\in  \mathbb{Z}$ and  $N_{L/K}(\beta) \in K$, they really differ by a unit $u$ 
		of $\mathbb{Z}\left[\xi\right].$ So, there exists a unit $u$ 
		$\in$ $U\left(\mathbb{Z}\left[ \xi\ 
		\right]\right)$ such that $u\cdot p^{h_{p}}=$$N_{L/K}(\beta).$ 
		According to Theorem \ref{teo24}  the
		symbol algebra $A=\left( \frac{\alpha,u\cdot p^{h_{p}}}{K,\xi}\right) $ splits.
		\item $p$ is not a primitive root modulo $q$.\\  Let's denote by $f$ the order of $p$ modulo   $q$.
		\newline
		According to Theorem \ref{teo22} we have $p\mathbb{Z}\left[\xi \right]
		=P_{1}\cdot...\cdot P_{r},$ where $P_{1},...,P_{r}$ are prime ideals in $\mathbb{Z}\left[ \xi\right]$ and $r=\frac{\phi\left(q\right)}{f}$ ($\phi$ as usual denotes the Euler's function).
		Since $\alpha$ is a $q$-th power residue modulo $p,$ it follows that $\alpha $ is a $q$-th
		 power residue modulo $P_{1},\ldots,P_{r}.$ 
		From Theorem \ref{teo23} we get:
		\\
		$p\mathcal{O}_{L}=P_{1}\mathcal{O}_{L}\cdot P_{2}\mathcal{O}_{L}\cdot...\cdot P_{r}\mathcal{O}_{L}$=\\
		$=P_{11}\cdot P_{12}\cdot....\cdot P_{1q}\cdot P_{21}\cdot P_{22}\cdot...\cdot P_{2q}\cdot...\cdot P_{r1}\cdot P_{r2}\cdot...\cdot P_{rq}$
		\\
		where $P_{1j},$ $P_{2j} , \ldots,P_{rj}$  
		are prime ideals in $\mathcal{O}_{L},$ with $j$$\in$$\left\lbrace 1,2,...,q\right\rbrace.$ 
		Since the ideals $P_{ij}$ are the prime divisors of
		$p_{i}\mathcal{O}_{L}$  ($i$$\in$$\left\lbrace 1,2,...,r\right\rbrace$, 
		 $j$$\in$$\left\lbrace 1,2,...,q\right\rbrace$)
		it follows that these ideals are   conjugate to each other under the action of the Galois group, 
		so their classes have the same order $h_{p}.$ Hence
		\\ $\left( p\mathcal{O}_{L}\right) ^{h_{p}}
		=\left( P_{11}\cdot P_{21}\cdot...\cdot P_{r1}\right)
		^{h_{p}}   \cdot\left( P_{12}\cdot P_{22}\cdot...\cdot P_{r2}\right)^{h_{p}}\cdot...\cdot \\ \left( P_{1q}\cdot P_{2q}\cdot...\cdot P_{rq}\right)
		^{h_{p}}$
		\\ 
		so there is an element $\gamma$$\in \mathcal{O}_{L}$ such that $p^{h_{p}}\mathcal{O}_{L}=N_{L/K}%
		\left(\gamma\right)\mathcal{O}_{L} .$ It follows that there exists a unit $u$ $\in$ $U\left(\mathbb{Z}\left[ \xi %
		\right]\right)$ such that $u\cdot p^{h_{p}}=$$N_{L/K}(\gamma).$ From Theorem  
		\ref{teo24}  it follows that the symbol algebras $%
		A=\left( \frac{\alpha, u\cdot p^{h_{p}}}{K,\xi}\right) $ splits.
			\end{itemize}

\item
According to Lemma 2.6    there exists a unit $u$ $\in$ $U\left(\mathbb{Z}\left[\xi\right]\right)$
		such that the symbol algebra $A=\left( \frac{\alpha, u\cdot p^{q}}{K,\xi }\right) $%
		splits. From this and from the first part of the proof, we obtain that there exists a unit $u$ $\in$ $U\left(\mathbb{Z}\left[ \xi\right]\right)$
		such that the symbol algebra $A=\left( \frac{\alpha, u\cdot p^{gcd\left( h_{p},q\right) }}{K,\xi }\right) $
		splits.
	\end{enumerate}
\end{proof}

The proof of the next proposition   
follows very closely  the proof of the first case of Proposition \ref{prop35}.
\begin{proposition}
\label{prop36}  {Let} $q$  {be an odd prime positive
integer and let}  $K=\mathbb{Q}\left( \xi \right)$ where $\xi $ {is a primitive root of order} $q$  {of unity.} {Let} $\alpha \in
K^{\ast }$  
 {and let } $L$ denote the Kummer field $K\left( \sqrt[q]{\alpha }\right)$.  
  {Let} $\pi$$\in$$\mathbb{Z}\left[ \xi \right]$  {be a prime element in the ring} $\mathbb{Z}\left[ \xi \right]$  {such that} $\alpha $  {is a} $q-$th  {power residue modulo} $\pi.$  {Let} $\mathcal{O}_L$  {be the ring of integers of the field} $L.$
 {Let} $Cl\left(L\right)$  {be the ideal class group of the ring} $\mathcal{O}_L,$  {let }$h_{L}$  {be the class number of} $L$ {and let}
 $h_{\pi}$ be the
 order in the 
 class group $Cl\left(L\right)$  of the class of a prime ideal in $\mathcal{O}_L$ that
 divides $\pi\mathcal{O}_L.$ 
{Then}:\\
a) {there exists a unit} $u$ $\in$ $U\left(\mathbb{Z}\left[\xi %
\right]\right)$
 { such that the symbol algebra}\\ $A=\left( \frac{\alpha ,u\cdot \pi^{h_{\pi}}}{K,\xi}\right) $%
 {\ splits;}\\
b) there exists a unit $u$ $\in$ $U\left(\mathbb{Z}\left[\xi %
\right]\right)$
such that the symbol algebra\\ $A=\left( \frac{\alpha, u\cdot \pi^{gcd\left( h_{\pi},q\right) }}{K,\xi}\right) $ splits.
\end{proposition}
In Corollary \ref{cor34} we obtained a sufficient condition for a symbol algebra $\left( \frac{\alpha ,u\cdot p^{h_{p}}}{K,\epsilon}\right)$ to split
 over the cyclotomic field $K=\mathbb{Q}\left(\epsilon\right)$, where $\epsilon$ is a primitive cubic root  of unity.
Next, let's ask ourselves if this condition is also necessary for the symbol algebra $\left( \frac{\alpha ,u\cdot p^{h_{p}}}{K,\epsilon}\right)$ to split.
The answer is affirmative, according to the following result:
\begin{theorem}
\label{th37}
 {Let} $\epsilon $  {be a primitive cubic
root  of unity and let }$K=\mathbb{Q}\left( \epsilon
\right)$. {Let} $\alpha \in
K^{\ast },$ let $p \neq 3$  {be}  {a prime rational integer,}  
 {and let} $L$ denote the Kummer field $K\left( \sqrt[3]{\alpha }\right)$.
  {Let} $\mathcal{O}_L$  {be the ring of integers of the field} $L.$
 {Let} $Cl\left(L\right)$  {be the ideal class group of the ring} $\mathcal{O}_L,$  {let } $h_{L}$  {be the class number of} $L$  {and let} $h_{p}$  {be the order of a class of a prime ideal in} $\mathcal{O}_L$, {which divides} $p\mathcal{O}_L,$  {in the group} $Cl\left(L\right).$  {Then, there exists a unit} $u$ $\in$ $U\left(\mathbb{Z}\left[ \epsilon %
\right]\right)$
 { such that the symbol algebra} $A=\left( \frac{\alpha, u\cdot p^{h_{p}}}{K,\epsilon }\right) $%
 {\ splits if and only if} $\alpha $  {is a cubic residue modulo} $p.$\\
\end{theorem}
\begin{proof}
We have proved  the sufficiency in Proposition \ref{prop35}, hence all we need to do
is to  prove the necessity.
Let's assume that there exists a unit $u$ of   $\mathbb{Z}\left[ \epsilon %
\right]$ such that the symbol algebra $A=\left( \frac{\alpha, u\cdot p^{h_{p}}}{K,\epsilon }\right)$ 
splits, where $\alpha $ is not a cubic residue modulo $p$.
This will lead to a contradiction. 
Clearly   $p\equiv$$1$ (mod $3$) and $\alpha^{\frac{p-1}{3}}$$\not\equiv$$1$ (mod $3$).
Since $\mathcal{O}_K=\mathbb{Z}\left[\epsilon\right]$ is a principal ideal domain, from Theorem \ref{teo22} and Theorem \ref{teo23} we get
\begin{equation}
p\mathcal{O}_K
=p_{1}\mathcal{O}_K \cdot p_{2}\mathcal{O}_K ,
\end{equation}
 where $p_{1},p_{2}$ are prime elements from $\mathcal{O}_K$ which remain primes in $\mathcal{O}_L,$ so
\begin{equation}
p\mathcal{O}_L
=P_{1}P_{2},
\end{equation}
where $P_{1}=p_{1}\mathcal{O}_L, \ P_{2}=p_{2}\mathcal{O}_L$
 are principal prime ideals of $\mathcal{O}_L$.
  It follows that $h_{p}=1,$ so
\begin{equation}
\left(p\mathcal{O}_L\right)^{h_{p}}=    p^{h_{p}}    \mathcal{O}_L   =P_{1}P_{2},
\end{equation}

We   show now that $u \cdot p=u \cdot p^{h_{p}}$, considered as an element of $K$,
 can not be the norm of an element of $L$, for any unit $u$ of $ \mathbb{Z}\left[\epsilon \right] $.
From Theorem \ref{teo24} it will follow that   the symbol algebra $A=\left( \frac{\alpha, u\cdot p^{h_{p}}}{K,\epsilon }\right) $ 
does not split,
and therefore $\alpha $ must be a cubic residue modulo $p.$

If  $u \cdot p = N_{L/K}(\beta)$ for some $\beta \in L$, then by definition
$u \cdot p = \beta \sigma(\beta) \sigma^2(\beta)$, where $\sigma$  generates  
the Galois group $Gal(L/K)$.
But then $u \cdot p^{h_{p}}    \mathcal{O}_L   =  ( \beta    \mathcal{O}_L )(  \sigma(\beta)  \mathcal{O}_L )(
  \sigma^2(\beta)    \mathcal{O}_L)$, 
and this contradicts the decomposition of $p\mathcal{O}_L$ as the product of two prime ideals 
$P_{1}$ and $   P_{2}$ of $\mathcal{O}_L$.
\end{proof}
In the paper \cite{Sa16} we obtained the following result:
\begin{proposition}
\label{prop38} \cite[Thm. 3.7]{Sa16}  {Let} $p$  {and} $q$  {be}  {%
prime positive integers such that} $p\equiv 1$ ({mod} $q$), let $\xi $ 
 {be a primitive root of order} $q$  {of unity and let }$K=%
\mathbb{Q}\left( \xi \right) $. {%
Then there is an integer} $\alpha $  {not divisible by} $p$ 
 {whose residue class mod} $p$  {does not belong to} $\left( 
\mathbb{F}_{p}^{\ast }\right) ^{q}$,  {and for every such an} $\alpha $%
,  {we have:}
\begin{itemize}
\item
 {the algebra} $A\otimes _{K}\mathbb{Q}_{p}$  {is a division
algebra over} $\mathbb{Q}_{p},$  {where} $A$  {is the symbol
algebra} $A=\left( \frac{\alpha ,p}{K,\xi }\right) ;$\newline
\item  {the symbol algebra} $A$  {is a division algebra over} $%
K.$
\end{itemize}
\end{proposition}
We consider now a Kummer field with class number $1$, for example   $L=\mathbb{Q}\left(\epsilon,\sqrt[3]{5}\right)$, 
where $\epsilon^{3}=1,  \  \epsilon\neq 1$. 
We consider the prime integers $17$ and $19.$ Now, $5$ is a cubic residue modulo $17,$ but
$5$ is not a cubic residue modulo $19.$ By using again  the computer algebra system MAGMA, we get that the norm equation $17=N_{L/\mathbb{Q}\left(\epsilon\right)}\left( a\right)$ has solutions, but the norm equation $N_{L/\mathbb{Q}\left(\epsilon\right)}\left( a\right)=19$
does not have a solution.
\smallskip\\

From  Proposition \ref{prop32} and Proposition \ref{prop38} we   obtain  in a very particular situation, i.e. when $L$ is a 
Kummer field of class number $1$, a necessary and sufficient condition for a symbol algebra to split over the third cyclotomic field.
\begin{proposition}
\label{prop39}  {Let} $\epsilon$  {be a primitive cubic root  of
unity and let }$K=\mathbb{Q}\left( \epsilon \right) $.
 {Let} $\alpha \in K^{\ast },$ let $p\neq 3$  {be a prime rational
integer,} {and let }$L=K\left( \sqrt[3]{\alpha }\right) $ 
 {be}  {a Kummer field with} $h_{L}=1$.   {Then, there exists a unit} $u$ $\in$ $U\left(\mathbb{Z}\left[ \epsilon\right]\right)$  { such that the symbol algebra} $A=\left( \frac{\alpha, u\cdot p}{K,\epsilon }\right) $%
 {splits if and only if} $\alpha $  {%
is a cubic residue modulo} $p.$
\end{proposition}
\begin{proof}
In order to prove the necessity, we note that, according to 
Remark \ref{rem25}, with $u=1$, the symbol algebra $A=\left( \frac{\alpha ,p}{K,\xi }\right)$ splits if and only if $A$ is not division algebra over $K.$ From Proposition 
\ref{prop38} it follows that $\alpha $ 
is a cubic residue modulo $p$ or $p$ is not congruent to $1$ modulo $3.$
But, if $p$ is not congruent to $1$ modulo $3$ and $p\neq 3,$ it follows that $p$$\equiv$$2$ (mod $3$) and 
this implies that $\left(\frac{\alpha}{p}\right)_{3}$ is equal to $1.$
Hence, from our previous results, we obtain that $\alpha $ is a cubic residue modulo $p.$

In order to prove the sufficiency,  we note that, if $\alpha $ is a cubic residue modulo $p,$ by
 applying Proposition \ref{prop32} with $h_L=1$  it follows that 
there exists a unit $u$ $\in$ $U\left(\mathbb{Z}\left[ \epsilon\right]\right)$  such that the symbol algebra $A=\left( \frac{\alpha, u\cdot p}{K,\epsilon }\right) $
splits.\\

Another way to prove Proposition \ref{prop39} is by using Theorem \ref{th37}, namely: if  $h_L=1,$ from  
 $h_p$ $\vert$
$h_L$ it follows that $h_p=1$, and by applying Theorem \ref{th37} we obtain that 
there exists a unit $u$ $\in$ $U\left(\mathbb{Z}\left[ \epsilon\right]\right)$   such that the symbol algebra $A=\left( \frac{\alpha, u\cdot p}{K,\epsilon }\right)$
splits if and only if $\alpha $ is a cubic residue modulo $p.$

\end{proof}
In the future   we will try to generalize the results contained in  
Theorem \ref{th37} and Proposition \ref{prop39} to Kummer fields
$L=\mathbb{Q}\left(\xi,\sqrt[l]{\alpha }\right)$,
where $l\geq 5$ is a prime integer and $\xi$ is a primitive root of order $l$ of the unity.
\section*{Acknowledgements}
Since  part of this work   has been done when the first author visited the University "G. D'Annunzio" of Chieti-Pescara, 
she wants to thank the  Department of Economic Studies of the University   for the hospitality and the support. In addition, she wants to thank Professor Claus Fieker
for the fruitful discussions about the computer algebra system MAGMA.\\
The authors thank anonymous referees for their comments and suggestions which helped us to improve this paper.

\end{document}